\documentclass[12pt,leqno,fleqn]{article}
\usepackage{amssymb, epsfig, amsmath, amsthm}
\usepackage{mathrsfs}    
\usepackage{color}       

\newcommand{\R}{\mathbb{R}}

\newcommand{\bp}{\boldsymbol \Phi}
\newcommand{\bpp}{\boldsymbol \Psi}

\newcommand{\bPhi}{\boldsymbol \Phi}

\newcommand{\supp}{\operatorname*{supp}}

\newcommand{\be}{\begin{equation}}
\newcommand{\ee}{\end{equation}}
\newcommand{\bea}{\begin{eqnarray}}
\newcommand{\eea}{\end{eqnarray}}
\newcommand{\bean}{\begin{eqnarray*}}
\newcommand{\eean}{\end{eqnarray*}}

\newcommand{\intl}{\int\limits}

\newcommand{\Beweisende}{\rule{0.2cm}{0.2cm}}

\newcommand{\intmw}{{\int\hspace{-830000sp}-\!\!}}

\newcounter{secnum}

\newtheorem{thm}{Theorem}[section]
\newtheorem{cor}[thm]{Corollary}

\theoremstyle{definition}

\title{On Liouville type theorems for the stationary \\
MHD and Hall-MHD systems} 
 
\author{Dongho Chae$^*$  and J\"{o}rg Wolf $^\dagger$\\
\ \\
 Department of Mathematics\\
Chung-Ang University\\
 Seoul 156-756, Republic of Korea\\
($*$)e-mail: dchae@cau.ac.kr \\
($^\dagger$)e-mail: jwolf2603@cau.ac.kr}
\date{}
\begin{document}

\maketitle
\begin{abstract}
In this paper we prove a Liouville type theorem for the stationary magnetohydrodynamics(MHD)  system in $\Bbb R^3$.
Let  $(v, B, p)$ be  a smooth solution to the  stationary MHD equations in $\Bbb R^3$. 
We show that if  there exist smooth matrix valued potential functions $\bp$, $\bpp$ such that  $ \nabla \cdot \bp =v$ and $\nabla \cdot \bpp= B$, 
whose $L^6$ mean oscillations have certain growth condition near infinity, namely 
$$ 
  \intmw_{B(r)} |\bp- \bp_{ B(r)} |^6 dx +\intmw_{B(r)} |\bpp- \bpp_{ B(r)} |^6 dx\le C r\quad  \forall 1< r< +\infty,
 $$
then $v=B= 0$ and $p=$constant.    With additional assumption of 
\[
r^{ -8}\intl_{B(r)} |B- B_{ B(r)}|^6 dx \rightarrow 0\quad   \text{as}\quad r \to +\infty,
\]
 similar result holds also  for the Hall-MHD system.
\ \\

  \noindent{\bf
keywords:} stationary magnetohydrodynamics equations, Liouville type theorem 
\end{abstract}
\section{Introduction}
\label{sec:-1}
\setcounter{secnum}{\value{section} \setcounter{equation}{0}
\renewcommand{\theequation}{\mbox{\arabic{secnum}.\arabic{equation}}}}

We consider the stationary magnetohydrodynamics equations in $ \R^{3}$.
$$
(MHD) \left\{ \aligned 
&-\Delta v + (v \cdot \nabla) v = - \nabla p + (B\cdot \nabla) B\quad  \text{in}\quad  \R^{3},\\
 & -\Delta B +(v\cdot \nabla) B -(B\cdot \nabla) v=0,\\
& \qquad \nabla \cdot v=\nabla \cdot B=0.
\endaligned
\right.
 $$

The system (MHD) describe the  steady  state of the physical system of fluid flows of plasma.
Here $v=v(x)=(v_1, v_2, v_3)$ is the velocity field of the fluid flows, $B=B(x)=(B_1, B_2, B_3)$ is the magnetic field, and $p=p(x)$
is the pressure of the  flows.  Note that if $B=0$, then the system (MHD)  reduces to the usual stationary Navier-Stokes system.
In this paper we study the Liouville type problem for the system (MHD).
The study is motivated by the similar Liouville problem for the stationary Navier-stokes equations, which is an active research area in the 
community of mathematical fluid mechanics(see e.g. \cite{gal, gil, koc, ser1, ser2, ser3, koz, pil, cham, cha1, cha2, cha3} and references therein). \\

We say  $\bp \in L^1_{loc} (\R^{n}; \R^{n\times  n}  )$ is  a potential function for  vector field $u\in L^1_{loc}  (\Bbb R^n)$ if  $ \nabla \cdot \bp=u$, where the derivative is in the sense of distribution. 
In \cite{ser1}  Seregin proved   Liouville type theorems for the Navier-Stokes equations under  assumption on the  potential functions of  the velocity filed, namely he showed the solution becomes trivial if the potential function of the velocity field belongs to $ BMO(\Bbb R^3)$ and  the velocity itself belongs to $L^6 (\Bbb R^3)$.
This was improved later in \cite{ser2},  removing the assumption for velocity  belonging to $L^6 (\Bbb R^3)$.
In recent paper \cite{sch}    Schulz    generalized Seregin's earlier result of \cite{ser1} to the the system(MHD), where  he proved that if the potential functions of the velocity and the magnetic field 
 belongs to $ BMO(\Bbb R^3)$,  while $(v,B) \in L^6 (\Bbb R^3)$,  then $v=B=0$. 
 In this paper we prove a Liouville theorem for the system  (MHD) under  more relaxed conditions  than \cite{sch}. 
 We allow some order of growth at infinity for the  mean oscillation for potentials of $v$ and $B$ without any integrability condition for them to obtain triviality of solutions.
 This result, on the other hand,  could be regarded as a generalization of the authors' previous result of \cite{cw} for the Navier-Stokes equations.\\

  For a measurable set $E \subset \Bbb R^n$  we denote by $|E |$ the $ n$-dimensional Lebesgue measure of $E$, and for $f\in L^1 (E)$ we use the notation
$$ f_E :=\intmw_{E}  f dx:= \frac{1}{|E|} \int_{E} fdx.
$$
Our aim in this paper is to prove the following:
\begin{thm}
 \label{thm1.1}
 Let $ (v, B, p) \in C^\infty(\R^{3} )\times  C^\infty(\R^{3} )\times C^\infty(\R^{3} )$ be a solution of (MHD). Suppose there exist $ \bp , \bpp\in  C^\infty(\R^{3}; \R^{3\times  3}  )$ such that 
 $ \nabla \cdot \bp = v$, $ \nabla \cdot \bpp = B$, and 
 \begin{equation}
  \intmw_{B(r)} |\bp- \bp_{ B(r)} |^6 dx +\intmw_{B(r)} |\bpp- \bpp_{ B(r)} |^6 dx\le C r\quad  \forall 1< r< +\infty.
 \label{1.2}
  \end{equation} 
 Then,  $ u =B= 0 $.
\end{thm}
If  $v, B\in BMO^{-1} (\Bbb R^3)$,  then there there exist  $\bp, \bpp$ such that  $ \nabla \cdot \bp = v$, $ \nabla \cdot \bpp = B$ and 
$$
 \sup_{r>1} \bigg( \intmw_{B(r)} |\bp- \bp_{ B(r)} |^6 dx +\intmw_{B(r)} |\bpp- \bpp_{ B(r)} |^6 dx\bigg) < +\infty, 
$$ 
and the condition \eqref{1.2} holds obviously, and therefore we have the following immediate corollary of the above theorem.
\begin{cor}
Let $ (v, B, p)$ be a smooth solution of  (MHD)  such that $v, B \in BMO^{-1} (\Bbb R^3)$, then   $ u =B= 0 $.
\end{cor}

 Next we consider the Hall-MHD system.
 $$
(HMHD)
\left\{ \aligned 
 &-\Delta u + (u \cdot \nabla) u = - \nabla p + (B\cdot \nabla) B\quad  \text{in}\quad  \R^{3},\\
 & -\Delta B +(u\cdot \nabla) B -(B\cdot \nabla) u=  \nabla \times ((\nabla \times B) \times B) ,\\
&\qquad \nabla \cdot u=\nabla \cdot B=0.
\endaligned
\right.
$$
This equations govern the the dynamics plasma flows of strong shear magnetic fields as in the solar flares, and have many important applications in the astrophysics.
We refer \cite{cw2}  and references therein for a recent mathematical approaches  and physical backgrounds for the system.  The following is our Liouville type theorem for (HMHD).
\begin{thm}
 \label{thm1.2}
 Let $ (v, B, p)$ be a smooth solution of (HMHD). 
Let us assume
 \be\label{1.3}
r^{ -8}\intl_{B(r)} |B- B_{ B(r)}|^6 dx \rightarrow 0\quad   \text{as}\quad r \to +\infty,
 \ee
and  there exist $ \bp , \bpp\in  C^\infty(\R^{3}; \R^{3\times  3}  )$ such that \eqref{1.2} holds. 
 Then,  $ u =B= 0 $.
\end{thm}

We believe that the additional condition  \eqref{1.3} in the case of (HMHD)  system  is sharp,  and it requires extra efforts
to deduce the Liouville property in Theorem\,\ref{thm1.2}. In fact, the iteration argument used below for the proof  is not standard, and would be of interest 
in itself.

 \section{Proof of Theorem 1.1 and Theorem 1.3}
 \label{sec:-2}
 \setcounter{secnum}{\value{section} \setcounter{equation}{0}
 \renewcommand{\theequation}{\mbox{\arabic{secnum}.\arabic{equation}}}} 

 In this section we prove  Theorem 1.1 and Theorem 1.3 in a unified fashion.
 Let us consider the following unified system with a constant $\alpha \ge 0$.
 \begin{align}\label{21}
 &-\Delta u + (u \cdot \nabla) u = - \nabla p + (B\cdot \nabla) B\quad  \text{in}\quad  \R^{3},\\
 \label{22} & -\Delta B +(u\cdot \nabla) B -(B\cdot \nabla) u= \alpha  \nabla \times ((\nabla \times B) \times B) ,\\
&\label{23} \qquad\qquad \nabla \cdot u=\nabla \cdot B=0.
 \end{align} 
 When $\alpha =0$ it reduces to (MHD), while for $\alpha=1$ it reduces to (HMHD). \\
 {\em We will use the assumption \eqref{1.3} only for the estimate of the term involving $\alpha$.}

\vspace{0.1cm}
 First let us show that  \eqref{1.2} implies the following estimate for the mean value. 
Let $ 1 <R< +\infty$. Let $ \varphi \in  C^\infty_c(B(R))$ denote a cut off function such that $ 0 \le \varphi  \le 1$ in $ B(R)$, 
 $ |\nabla \varphi | \le c R^{ -1}$ and 
 \begin{equation}
\intl_{B(R)} \varphi dx \ge c R^3.   
 \label{2.0a}
  \end{equation} 
 For a function $ f\in L^1(B(R))$ we define the corresponding  mean value 
\[
f_{ B(R), \varphi } = \frac{1}{\intl_{B(R)} \varphi  dx } \intl_{B(R)} f\varphi  dx. 
\]
Using H\"older's inequality and observing  \eqref{2.0a}, we get for all $ f\in L^p(B(R)), 1 \le p \le +\infty$ the estimate 
\begin{equation}
|f_{ B(R), \varphi }| \le c\|f\|_{ L^p(E)} \quad   \text{for every measurable $B(R)\supset E  \supset \supp(\varphi )$}.  
\label{2.0b}
 \end{equation} 
On the other hand, for every  constant function $ a\in \R$ it holds  $ a_{ B(R), \varphi } =a$.  This provides us with  
\[
f - f_{ B(R), \varphi } = f-a - (f-a)_{ B(R), \varphi }.
\]
Thus, applying  \eqref{2.0b}, we get for all $ a\in \R$ and for every measurable $ B(R)\supset E  \supset \supp(\varphi )$, 
\begin{equation}
\| f - f_{ B(R), \varphi }\|_{ L^p(E)} \le c\| f- a\|_{ L^p(E)}. 
\label{2.0c}
 \end{equation} 
 In particular,  \eqref{2.0c} with $  a= f_{ E}$ gives 
\begin{equation}
\| f - f_{ B(R), \varphi }\|_{ L^p(E)} \le c\| f- f_{ E}\|_{ L^p(E)}. 
\label{2.0d}
 \end{equation} 
 This, inequality will be used below for the application of both, the Poincar\'e inequality and the Sobolev-Poincar\'e
inequality.  

Using  the generalized mean we have introduced above,   we now get from  \eqref{1.2} the following  estimate for $ u$ and $ B$ 
respectively, 
 \begin{equation}
|u_{ B(R), \varphi }| \le c R^{ - \frac{5}{6}},\quad  |B_{ B(R), \varphi }| \le c R^{ - \frac{5}{6}}. 
 \label{2.0e}
  \end{equation} 
 Indeed, applying integration by parts, Jensen's inequality, and observing  \eqref{1.2} along with  \eqref{2.0a}, we find 
 \begin{align*}
   |u_{ B(R), \varphi }| &\le c R^{ -3}  \bigg|\intl_{B(R)} \nabla \cdot  (\bPhi- (\bPhi)_{ B(R)}) \varphi dx \bigg|
 \\
 &= c\bigg|\intmw_{B(R)} (\bPhi- \bPhi_{ B(R)}) \cdot \nabla \varphi dx \bigg|
 \le R^{ -1} \bigg(\intmw_{B(R)} |\bPhi- \bPhi_{ B(R)}|^6 dx \bigg)^{ \frac{1}{6}}
 \\
 &\le c R^{ - \frac{5}{6}}. 
 \end{align*}
 Whence,  \eqref{2.0e} for $ u$. The estimate for $ B$ follows by the same reasoning from  \eqref{1.2}.

\vspace{0.1cm}
For the sake of notational simplicity in our discussion below we use the following abbreviations
 \begin{align*}
\Theta  ( r) &:= \alpha \bigg(r^{ -8} \intl_{B(r)} |B- B_{ B(r)}|^6 dx\bigg)^{ \frac{1}{6}},    
\\[0.2cm]
G(r) &:= \intl_{B(r)} |\nabla u|^2 dx + \intl_{B(r)} |\nabla B|^2 dx,  \quad 0< r< +\infty.
\end{align*}  
Note that in case of $ \alpha =1$, in view of  condition \eqref{1.3} for $ B$ it holds 
\begin{equation}
 \lim_{r \to \infty}  \Theta  (r) = 0. 
\label{2.22a}
 \end{equation} 
Otherwise in case $ \alpha =0$ it holds $ \Theta \equiv 0$.  

Furthermore, during the proof below we make frequently use of the following elementary estimate, for all  $ \alpha , \beta, \gamma  \in \R$ with 
$ 0 \le \alpha  \le \beta \le \gamma $ and for all $ 1 \le \rho < R < +\infty$
\begin{equation}
 R^{ \alpha } (R-\rho )^{ -\beta } \le R^{ \gamma } (R-\rho)^{ -\gamma }.
\label{2.0f}
 \end{equation} 
 Indeed, since $ \alpha \le \beta $ we see that $ R^{ \alpha } \le R^{ \beta }$ for all $ R \ge 1$.  On the other hand,  by means of 
 $ \beta \le \gamma $  and $ R (R-\rho )^{ -1} >1$ we infer that $ R^{ \beta } (R-\rho )^{ -\beta } \le R^{ \gamma  } (R-\rho )^{ -\gamma  }  $
 Accordingly, for all  $ 1 \le R < \rho < +\infty$ we have
 \begin{align*}
  R^{ \alpha } (R-\rho )^{ -\beta }  &\le    R^{ \beta } (R-\rho )^{ -\beta }    \le  R^{ \gamma  } (R-\rho )^{ -\gamma  }.    
  \end{align*}

Let   $ 1 < r< +\infty$ be arbitrarily chosen, but fixed. Let   $ r \le  \rho < R \le 2r $. We set $ \overline{R}= \frac{R+ \rho }{2}$.  Let $ \zeta \in  C^\infty(\R^{n} )$ be a cut off function, which is radially non-increasing   with $ \zeta =1$ on $B(\rho )$ and 
$\zeta=0$ on $ \Bbb R^3\setminus B(\overline{R})$ satisfying 
$ |\nabla \zeta | \le c (R-\rho )^{ -1}$ and $ |D ^2 \zeta | \le c (R-\rho )^{ -2}$. We multiply  \eqref{21} by $ u \zeta ^2$,  integrate it over $ B(R)$,  and  then 
we  multiply  \eqref{22} by $ B \zeta ^2$,  integrating it over $ B(R)$, and add them together. 
Then, after  integration by parts we obtain
\begin{align}
\hspace*{-2cm} G(\rho ) &= \intl_{B(\rho )} |\nabla u|^2 dx + \intl_{B(\rho )} |\nabla B|^2 dx  \le  \frac{1}{2}\intl_{B(\overline{R})} |u|^2 \Delta \zeta ^2 dx \cr
 &+ \frac{1}{2}\intl_{B(\overline{R})} |u|^2 u \cdot \nabla \zeta ^2  dx+ \intl_{B(\overline{R})} (p - p_{ B(\overline{R})}) u \cdot \nabla \zeta ^2 dx  - \intl_{B(\overline{R})}  B\cdot u  ( B\cdot \nabla ) \zeta^2 dx \cr
& +  \frac{1}{2}\intl_{B(\overline{R})} |B|^2 \Delta \zeta ^2 dx   + \frac{1}{2}\intl_{B(\overline{R})} |B|^2 u \cdot \nabla \zeta ^2  dx    + \intl_{B(\overline{R})}  u\cdot B   ( B\cdot \nabla ) \zeta^2 dx \cr
& -\alpha \intl_{B(\overline{R})} ((\nabla \times B) \times B) \cdot B\times \nabla \zeta^2 dx 
\cr
&\le c (R- \rho )^{ -2} \intl_{B(\overline{R})} |u|^2  dx+  c (R- \rho )^{ -1} \intl_{B(\overline{R})} |u|^3dx  
\cr
& + c   (R- \rho )^{ -1} \intl_{B(\overline{R})} |p - p_{ B(\overline{R}) }| |u| dx+ c   (R- \rho )^{ -1} \intl_{B(\overline{R})} |B|^2 |u| dx  
\cr
 &+ c (R- \rho )^{ -2} \intl_{B(\overline{R})} |B|^2  dx +\alpha  c   (R- \rho )^{ -1}    \intl_{B(\overline{R})} |B|^2 |\nabla B| dx\cr
  &= I+II+ III +IV+V+VI. 
 \label{2.1}
\end{align}

In order to estimate $ I$ we choose another   cut off function  $ \psi \in  C^\infty(\R^{3} )$, 
which is radially non-increasing   with $ \psi =1$ on $ B(\overline{R})$ and 
$\psi=0$ on $ \Bbb R^3\setminus B(R)$ satisfying 
  $ |\nabla \psi  | \le c (R-\rho )^{ -1}$.
 Recalling that $ u = \nabla \cdot \bp$, applying integration by parts   and Cauchy-Schwarz' inequality,  
we find 
\begin{align*}
\intl_{B(R)} |u|^2\psi ^2 dx &= \intl_{B(R)} \partial _i (\Phi_{ ij}- (\Phi_{ ij})_{ B(R)}) u_{ j} \psi^2 dx  
\\
&  =  -\intl_{B(R)} (\Phi_{ ij}- (\Phi_{ ij})_{ B(R)}) \partial _i u_{ j} \psi^2 dx -  \intl_{B(R)} (\Phi_{ ij}- (\Phi_{ ij})_{ B(R)}) u_{ j} \partial _i\psi^2 dx\\
& \le c\bigg(\intl_{B(R)}| \bp - \bp_{ B(R)}|^2 dx \bigg)^{  \frac{1}{2}} \bigg(\intl_{B(R)}| \nabla u|^2 dx \bigg)^{  \frac{1}{2}}
\\
&\qquad + c (R- \rho )^{ -1}\bigg(\intl_{B(R)}| \bp - \bp_{ B(R)}|^2 dx \bigg)^{  \frac{1}{2}} \bigg(\intl_{B(R)}| u|^2\psi ^2 dx \bigg)^{  \frac{1}{2}}.
\end{align*}   
 Using Jensen's inequality and Young's inequality, we obtain 
 \begin{align*}
   \intl_{B(R)} |u|^2\psi ^2 dx & \le c R\bigg(\intl_{B(R)}| \bp - \bp_{ B(R)}|^6 dx \bigg)^{  \frac{1}{6}} \bigg(\intl_{B(R)}| \nabla u|^2 dx \bigg)^{  \frac{1}{2}}
\\
&\qquad + c R^2  (R- \rho )^{ -2}\bigg(\intl_{B(R)}| \bp - \bp_{ B(R)}|^6 dx \bigg)^{  \frac{1}{3}}.
 \end{align*}  
 Observing  \eqref{1.2},  and using Young's inequality together with  \eqref{2.0f} with $ \gamma =11$, we infer 
 \begin{align*}
 I &\le  c R^{ \frac{5}{3}} (R- \rho )^{ -2} \bigg(\intl_{B(R)}| \nabla u|^2 dx \bigg)^{  \frac{1}{2}} + cR^{ \frac{10}{3}} (R-\rho )^{ -4}
 \\
 &   \le \frac{1}{100} G(R) + cR^{11} (R-\rho )^{ -11}. 
 \end{align*}  
 Similarly,
 \be
 V \le \frac{1}{100} G(R) + cR^{11} (R-\rho )^{ -11}.
 \ee
In order to estimate $II$ we  first estimate the $ L^3$ norms of $ u$   and $w$ as follows.
\begin{align*}
&\intl_{B(R)} |u|^3\psi ^3 dx = \intl_{B(R)} \partial _i (\Phi_{ ij}- (\Phi_{ ij})_{ B(R)}) u_{ j} |u|\psi^3 dx  
\\
&  = - \intl_{B(R)} (\Phi_{ ij}- (\Phi_{ ij})_{ B(R)}) \partial _i (u_{ j} |u|)\psi^3 dx -    \intl_{B(R)} (\Phi_{ ij}- (\Phi_{ ij})_{ B(R)}) u_{ j}|u| \partial _i\psi^3 dx
\end{align*}
\begin{align*}
& \le c\bigg(\intl_{B(R)}| \bp - \bp_{ B(R)}|^6 dx \bigg)^{  \frac{1}{6}} \bigg(\intl_{B(R)}| u|^3 \psi ^3 dx \bigg)^{  \frac{1}{3}}\bigg(\intl_{B(R)}| \nabla u|^2 dx \bigg)^{  \frac{1}{2}}
\\
&\qquad + c (R- \rho )^{ -1}\bigg(\intl_{B(R)}| \bp - \bp_{ B(R)}|^3 dx \bigg)^{  \frac{1}{3}} \bigg(\intl_{B(R)}| u|^3\psi ^3 dx \bigg)^{  \frac{2}{3}}.
\end{align*}   
Using  Young's inequality, we get 
\begin{align}
\intl_{B(R)} |u|^3\psi ^3 dx 
& \le c\bigg(\intl_{B(R)}| \bp - \bp_{ B(R)}|^6 dx \bigg)^{  \frac{1}{4}}\bigg(\intl_{B(R)}| \nabla u|^2 dx \bigg)^{  \frac{3}{4}}
\cr
&\qquad + c R^{ \frac{3}{2}}(R- \rho )^{ -3}\bigg(\intl_{B(R)}| \bp - \bp_{ B(R)}|^6 dx \bigg)^{  \frac{1}{2}}. 
 \label{2.3a}
\end{align}  
Multiplying  \eqref{2.3a} by $ (R-\rho )^{ -1}$ combined with the hypothesis \eqref{1.2}, and using  \eqref{2.0f} $ \gamma =4$, we infer 
\begin{equation}
 (R- \rho )^{ -1}\intl_{B(R)} |u|^3\psi ^3 dx \le R(R- \rho )^{ -1} G(R)^{ \frac{3}{4}}+ cR^4(R- \rho )^{ -4}. 
\label{2.3b}
 \end{equation}

Applying Young's inequality,  and  \eqref{2.0f} with $ \gamma =11$ one has 
\begin{align*} 
II \le  \frac{1}{100} G(R)+c R^{11} (R-\rho)^{-11}.
\end{align*} 
Similarly to \eqref{2.3a},  we also obtain 
\begin{align}
 (R- \rho )^{ -1}\intl_{B(R)} |B|^3\psi ^3 dx \le R(R- \rho )^{ -1} G(R)^{ \frac{3}{4}}+ cR^4(R- \rho )^{ -4}. 
 \label{2.3c}
\end{align}   
Next, using \eqref{2.3a}, \eqref{2.3b},  we shall estimate $\intl_{B(R)} |u| |B|^2\psi ^3 dx$. Applying Young's inequalities, and \eqref{2.3a}, \eqref{2.3b}
we have
\begin{align*}
(R- \rho )^{ -1}\intl_{B(R)} |u| |B|^2 \psi ^3 dx & \le (R- \rho )^{ -1} \int_{B(R)}( |u|^3 + |B|^3) \psi^3 dx
\\
&\le cR(R- \rho )^{ -1} G(R)^{ \frac{3}{4}}+c R^{4} (R-\rho)^{-4}.
 \end{align*}

Applying Young's inequality, and again using  \eqref{2.0f} with $ \gamma =11$,  we arrive at 
 \begin{align*}
IV \le     \frac{1}{100}G(R)+cR^{ 11} (R-\rho )^{ -11}. 
 \end{align*}  
We now  estimate $III$. Using H\"older's inequality and Young's inequality,  we infer 
\begin{align}
III&\le c (R- \rho )^{ -1} \intl_{B(\overline{R})} |p- p_{ B(\overline{R})}|^{ \frac{3}{2}} dx   + 
c (R- \rho )^{ -1} \intl_{B(\overline{R})} |u|^{3} dx.   
 \label{2.20} 
\end{align}

We are now going to estimate the pressure term.

 For this purpose let us define the   functional   $ F\in W^{-1,\, \frac{3}{2}}(B(\overline{R})) $, by means of 
\[
 \langle F, \varphi \rangle =  \intl_{B(\overline{R})} (\nabla u - u \otimes u +B\otimes B) :\nabla \varphi dx,\quad  \varphi \in W^{1,\,3}_0(B(\overline{R})).  
\]
Since $ (u, w, p)$ is a solution to  \eqref{21}, it follows that  
\[
 \langle F, \varphi \rangle = \intl_{B(\overline{R})} (p- p_{ B(\overline{R})}) \nabla \cdot \varphi dx \quad \forall \varphi \in W^{1,\,3}_0(B(\overline{R})).
\]
Thus, we have $ F = - \nabla (p- p_{ B(\overline{R})})$ with $ p-p_{ B(\overline{R})}\in L^{ \frac{3}{2}}(B(\overline{R}))$ and $ \intl_{B(\overline{R})}
(p- p_{ B(\overline{R})} )dx =0$. Consulting Sohr \cite[Lemm\,2.1.1]{sohr}, we get the estimate 
 \begin{align}
 \label{2.21}\intl_{B(\overline{R})} |p- p_{ B(\overline{R})}|^{ \frac{3}{2}} dx   &\le c \|F\|^{ \frac{3}{2}}_{ W^{-1,\, \frac{3}{2}} (B(\overline{R}))},
\end{align}
with a constant $ c>0$ independent of $ \overline{R}$. On the other hand, we estimate by the aid of H\"older's inequality 
\begin{align}
\|F\|^{ \frac{3}{2}}_{ W^{-1,\, \frac{3}{2}} (B(\overline{R}))}
&\le \|\nabla u - u \otimes u + B\otimes B)\|^{ \frac{3}{2}}_{ L^{ \frac{3}{2}}(B(\overline{R}))}
\cr
&\le c  R^{ \frac{3}{4}} \bigg(\intl_{B( \overline{R})} |\nabla u|^2 dx \bigg)^{ \frac{3}{4}} + c  \intl_{B( \overline{R})} |u|^3 dx+ c  \intl_{B( \overline{R})} |B|^3 dx.
 \label{2.21a}
\end{align}

Combining  \eqref{2.20},  \eqref{2.21}, and  \eqref{2.21a},  and using Young's inequality along with  \eqref{2.0f}  for $ \gamma =11$, we obtain 
\be\label{2.21c} 
III\le \frac{1}{100} G(R)  +c R^{ 11}(R- \rho )^{ -11}+ c  (R- \rho )^{ -1}\intl_{B(\overline{R})} |u|^3 dx+c  (R- \rho )^{ -1}\intl_{B(\overline{R})} |B|^3 dx.
\ee
Note that \eqref{2.3a}, \eqref{2.3b} provide us with the estimate
\begin{align}\label{2.21f} 
&  (R- \rho )^{ -1}  \intl_{B(R)} |u|^3 \psi ^3 dx+  (R-\rho )^{ -1} \intl_{B(R)} |B|^3 \psi ^3 dx
\cr
&\qquad  \le c R(R- \rho )^{ -1}  G(R)^{ \frac{3}{4}}
    +c  R^{4} (R-\rho )^{-4}.
   \end{align}
  Inserting the estimate of \eqref{2.21f} into  \eqref{2.21c}, applying Young's inequality,    and using  \eqref{2.0f}, we find 
 \begin{align*}
III&\le   \frac{1}{100} G(R)  +cR^{ 11} (R-\rho )^{ -11}. 
  \end{align*}

  It remains to estimate $ VI$. Applying H\"older's inequality, we estimate 
\begin{align*}
VI &\le c\alpha  (R-\rho )^{ -1} \intl_{B(R)} |B|^2 |\nabla B| \psi  dx 
\\
& \le c \alpha (R-\rho )^{ -1}\bigg(\intl_{B(R)} |B|^6 dx\bigg)^{ \frac{1}{6}} 
\bigg(\intl_{B(R)} |B|^3 \psi ^3 dx\bigg)^{ \frac{1}{3}}\bigg(\intl_{B(\overline{R})} |\nabla B|^2 dx\bigg)^{ \frac{1}{2}} 
\\
& \le c\alpha  (R-\rho )^{ -1}\bigg(\intl_{B(R)} |B - B_{ B(R), \psi }|^6 dx\bigg)^{ \frac{1}{6}} 
\bigg(\intl_{B(R)} |B|^3 \psi ^3 dx\bigg)^{ \frac{1}{3}}\bigg(\intl_{B(\overline{R})} |\nabla B|^2 dx\bigg)^{ \frac{1}{2}} 
\\
& \quad  + c \alpha R^{ \frac{1}{2}}(R-\rho )^{ -1} |B_{ B(R), \psi }| 
\bigg(\intl_{B(R)} |B|^3 \psi ^3 dx\bigg)^{ \frac{1}{3}}\bigg(\intl_{B(\overline{R})} |\nabla B|^2 dx\bigg)^{ \frac{1}{2}}.  
\end{align*}
Applying  \eqref{2.0d} with $ E= B(R)$, and  \eqref{2.0e} both with $ \varphi = \psi $,   we infer  from the above estimate 
\begin{align}
VI  & \le c \alpha  (R-\rho )^{ -1}\bigg(\intl_{B(R)} |B - B_{ B(R) }|^6 dx\bigg)^{ \frac{1}{6}} 
\bigg(\intl_{B(R)} |B|^3 \psi ^3 dx\bigg)^{ \frac{1}{3}}\bigg(\intl_{B(\overline{R})} |\nabla B|^2 dx\bigg)^{ \frac{1}{2}} 
\cr
& \quad  + c\alpha  R^{ -\frac{1}{3} }(R-\rho )^{ -1}  
\bigg(\intl_{B(R)} |B|^3 \psi ^3 dx\bigg)^{ \frac{1}{3}}\bigg(\intl_{B(\overline{R})} |\nabla B|^2 dx\bigg)^{ \frac{1}{2}}
\cr
& \le c R^{ \frac{4}{3}}(R-\rho )^{ - \frac{2}{3}} \Theta  (R) 
\bigg( (R-\rho )^{ -1}\intl_{B(R)} |B|^3 \psi ^3 dx\bigg)^{ \frac{1}{3}}\bigg(\intl_{B(\overline{R})} |\nabla B|^2 dx\bigg)^{ \frac{1}{2}} 
\cr
& \quad  + c R^{ -\frac{1}{3} }(R-\rho )^{ - \frac{2}{3}}  
\bigg((R-\rho )^{ -1}\intl_{B(R)} |B|^3 \psi ^3 dx\bigg)^{ \frac{1}{3}}\bigg(\intl_{B(\overline{R})} |\nabla B|^2 dx\bigg)^{ \frac{1}{2}}.
 \label{2.21d}
  \end{align}

Combining \eqref{2.21f} with \eqref{2.21d}, and applying Young's inequality,  and using  \eqref{2.0f} several times with $ \gamma =11$,  we find 
  \begin{align*}
  VI & \le  c R^{ \frac{4}{3}}(R-\rho )^{ - \frac{2}{3}} \Theta (R)  \Big( R^{ \frac{1}{3}}(R-\rho )^{ - \frac{1}{3}} G(R)^{ \frac{1}{4}} + R^{ \frac{4}{3}} (R-\rho )^{ - \frac{4}{3}}\Big) G(R)^{ \frac{1}{2}}
  \\
&\qquad   + c  R^{ -\frac{1}{3} }(R-\rho )^{ - \frac{2}{3}}  
\Big( R^{ \frac{1}{3}}(R-\rho )^{ - \frac{1}{3}} G(R)^{ \frac{1}{4}} + R^{ \frac{4}{3}} (R-\rho )^{ - \frac{4}{3}}\Big) G(R)^{ \frac{1}{2}}
\\
 &\le cR^{ \frac{20}{3}}(R-\rho )^{ - 4} \Theta^4 (R) + c R^{ \frac{16}{3}}(R-\rho )^{ - 4} \Theta^2 (R)
 +c (R- \rho )^{ -4} +c R^2 (R- \rho )^{ -4}
 \\
& \qquad + \frac{1}{100} G(R)
\\
& \le   \frac{1}{100} G(R) +cR^{ \frac{20}{3}}(R-\rho )^{ - 4} \Theta^4 (R) + cR^{ 11} (R-\rho )^{ -11}
  \end{align*}
    
Inserting the estimates of $ I,\cdots, VI$ into the right hand side of  \eqref{2.1}, again using $ R \ge 1$,  and applying Young's inequality,  we are led to  
  \begin{align}\label{iter1}
G(\rho )   &\le \frac{1}{2} G(R)  +  R^{ \frac{20}{3}}(R-\rho )^{ - 4} \Theta^4  (R)  + cR^{11} (R-\rho )^{ - 11}. 
  \end{align}

Applying the iteration Lemma in \cite[V.\,Lemma\,3.1]{gia},  from \eqref{iter1} we deduce that 
\be\label{2.4a}
G(\rho )\le cR^{ \frac{20}{3}}(R-\rho )^{ - 4} \Theta^4 (R) + cR^{ 11} (R-\rho )^{ - 11}, \quad \forall r\le \rho <R \le 2r.  
\ee
In the case  $ \alpha =0$ (MHD equations) we see that the first term on the right-hand side vanishes, and thus, $ G$ is bounded. In case  $ \alpha =1$ 
the inequality  \eqref{2.4a} with $ \rho =r$ and $ R=2r$  gives  
\begin{equation}
G(r) \le c r^{ \frac{8}{3}}\quad  \forall 1 \le r < +\infty. 
\label{2.4aa}
 \end{equation} 
Since our aim is to show that $ \nabla u$ and $ \nabla B$ are in $ L^2(\R^{3} )$  we   still need to improve the above   estimate. To do this we proceed as follows.  Let $ 0< \tau < 1$ be sufficiently small, specified below.  
Observing  \eqref{2.22a},  we may choose $ r_0>1$ such that $ \Theta (r) \le \tau ^2$ for all $ r \ge r_0$, 

Let $ r \ge r_0$ be arbitrarily chosen, but fixed. 
Taking into account   \eqref{2.4a}, we repeat the estimation of $ VI$.  Starting from the first inequality in  \eqref{2.21d} applying Sobolev's-Poincar\'e 
inequality, we find  
\begin{align}
VI  & \le c (R-\rho )^{ -1}\bigg(\intl_{B(R)} |B|^3 \psi ^3 dx\bigg)^{ \frac{1}{3}}\bigg(\intl_{B(R)} |\nabla B|^2 dx \bigg)^{ \frac{1}{2}}\bigg(\intl_{B(\overline{R})} |\nabla B|^2 dx\bigg)^{ \frac{1}{2}}
\cr
& \qquad  + c R^{ -\frac{1}{3} }(R-\rho )^{ -1}  
\bigg(\intl_{B(R)} |B|^3 \psi ^3 dx\bigg)^{ \frac{1}{3}}\bigg(\intl_{B(R)} |\nabla B|^2 dx\bigg)^{ \frac{1}{2}}
\cr
& \le c (R-\rho )^{ - \frac{2}{3}}\bigg( (R-\rho )^{ -1}\intl_{B(R)} |B|^3 \psi ^3 dx\bigg)^{ \frac{1}{3}} G(R)^{ \frac{1}{2}}G(\overline{R})^{ \frac{1}{2}} 
\cr
& \qquad  + c R^{ -\frac{1}{3} }(R-\rho )^{ - \frac{2}{3}}  
\bigg((R-\rho )^{ -1}\intl_{B(R)} |B|^3 \psi ^3 dx\bigg)^{ \frac{1}{3}} G(R)^{ \frac{1}{2}}.
 \label{2.21j}
  \end{align}
Inserting the estimate  \eqref{2.21f} into the right-hand side of  \eqref{2.21j}, we arrive at 
\begin{align*}
VI &\le  c (R-\rho )^{ - \frac{2}{3}}\Big( R^{ \frac{1}{3}}(R-\rho )^{ - \frac{1}{3}} G(R)^{ \frac{1}{4}} + R^{ \frac{4}{3}} (R-\rho )^{ - \frac{4}{3}}\Big) 
G(\overline{R})^{ \frac{1}{2}} G(R)^{ \frac{1}{2}}
\\
&\qquad + c R^{ - \frac{1}{3}} (R-\rho )^{ - \frac{2}{3}}\Big( R^{ \frac{1}{3}}(R-\rho )^{ - \frac{1}{3}} G(R)^{ \frac{1}{4}} + R^{ \frac{4}{3}} (R-\rho )^{ - \frac{4}{3}}\Big) G(R)^{ \frac{1}{2}}
\\
& =c R^{ \frac{1}{3}} (R-\rho )^{ - 1} G(\overline{R})^{ \frac{1}{4}} G(\overline{R})^{ \frac{1}{4}} G(R)^{ \frac{3}{4}}+ c R^{ \frac{4}{3}} (R-\rho )^{ - 2} G(\overline{R})^{ \frac{1}{4}} 
G(R)^{ \frac{3}{4}}
\\
& \qquad +c (R-\rho )^{ - 1} G(R)^{ \frac{3}{4}}+ c R (R-\rho )^{ - 2}  G(R)^{ \frac{1}{2}}. 
\end{align*}

 Estimating  $ G(\overline{R})^{ \frac{1}{4}}$ by means of  \eqref{2.4a} with $ \rho =\overline{R}$, using  \eqref{2.0f} with $ \gamma =11$, and applying Young's inequality, we are provided with 
\begin{align*}
VI &\le c R^{ \frac{5}{3}} (R- \rho )^{ - \frac{5}{3}}   \Theta  (R)  G(\overline{R})^{ \frac{1}{4}}G(R)^{ \frac{3}{4}}+ c R^{ \frac{5}{3}}(R-\rho )^{- \frac{7}{3}} G(\overline{R})^{ \frac{1}{4}} 
G(R)^{ \frac{3}{4}}. 
\\
&\qquad + cR^{ \frac{32}{3}} (R-\rho )^{ -\frac{32}{3}} + cR^{ \frac{4}{3}} (R-\rho )^{- \frac{10}{3}}
+ c R^{ 4} (R-\rho)^{ -4} + \frac{1}{100} G(R)
\\
& \le  \frac{1}{50} G(R) + c R^{ 11} (R- \rho )^{ - 11}   \Theta  (R)  G(R)+ c R^{ 3}(R-\rho )^{-3} G(R)^{ \frac{3}{4}} 
\\
&\qquad +cR^{ 11}(R-\rho )^{-11}
\\
&\le \frac{1}{25} G(R)+  c R^{ 11} (R- \rho )^{ - 11}   \Theta  (R)  G(R) +  cR^{11} (R-\rho )^{ -11}. 
  \end{align*}
Again inserting the estimates of $ I,\cdots, VI$ into the right hand side of  \eqref{2.1}, and applying Young's inequality,  we are led to  
  \begin{align}\label{iter2}
G(\rho )   &\le \frac{1}{2} G(R)  +  R^{ 11}(R-\rho )^{ - 11} \Theta  (R)G(R)  + cR^{11} (R-\rho )^{ -11}.
  \end{align}

Once more applying the iteration Lemma in \cite[V.\,Lemma\,3.1]{gia},  from \eqref{iter2} we deduce that 
\be\label{2.4d}
 G(\rho )\le  c R^{ 11} (R- \rho )^{ -11}  (\tau^2  G(2r) +  1), \quad \forall r\le \rho <R \le 2r.
\ee

In particular, this inequality with $ \rho = r$ and $ R= 2r$   \eqref{2.4d} reads
\be\label{2.4e}
G(r )\le  c  \tau^2  G(2r) + c.   
\ee
with an absolute constant $ c>0$. We take $ 0< \tau < 1$ such that $ c\tau \le 1 $, and iterate   \eqref{2.4e} $ k$-times starting with some $ R \ge r_0$. 
This gives 
\be\label{2.4e}
  G( R )\le    \tau^{ k}  G(2^k R) +  c  \sum_{i=0}^{k-1} \tau ^i.   
\ee
We may choose $ 0< \tau \le 2^{- 3}$.  With this choice along with  \eqref{2.4aa} we get 
\be\label{2.4e}
  G( R )\le    2^{-3 k}  G(2^k R) +  2c  \le   c2^{- 3 k}R^{ \frac{8}{3}} 2^{ \frac{8}{3} k } 
  + 2c\le R^{ \frac{8}{3}} 2^{ - \frac{1}{3} k} + 2c.  
 \ee
After letting $ k \rightarrow +\infty$  and then  passing $R\to +\infty$, we find 
\begin{equation}
\int_{\Bbb R^3} |\nabla u|^2 dx+  \int_{\Bbb R^3} |\nabla B|^2 dx  <+\infty .
\label{2.4}
 \end{equation} 
 Note that from \eqref{2.21f} with $ \rho =r$ and $ R= 2\rho $, using  \eqref{2.4},  provides us with  the estimate 
 \begin{equation}
 r^{ -1}  \bigg( \intl_{B(r)} |u|^3 dx  + \intl_{B(r)} |B|^3 dx \bigg)   \le c  \quad  \forall 1< r< +\infty. 
 \label{2.5}
  \end{equation} 

Next, we claim that 
\begin{equation}
 r^{ -1} \bigg( \intl_{B(3r) \setminus B(2r)} |u|^3 dx +  \intl_{B(3r) \setminus B(2r)} |B|^3 dx\bigg)  = o(1)\quad  \text{as}\quad r \rightarrow +\infty. 
\label{2.6}
 \end{equation} 
 Let $ \psi \in  C^\infty( \R^{3} )$ be a cut off function for the annulus $ B(3r) \setminus B(2r)$ in $ B(4r) \setminus B(r)$, 
i.e. $ 0 \le \psi \le 1$ in $ \R^{3} $, $ \psi = 0 $ in $ \R^{3} \setminus (B(4r) \setminus B(r))$, $ \psi =1$ on $ B(3r) \setminus B(2r)$ and  
 $ |\nabla \psi | \le c r^{-1}$  .  Recalling   that $ u = \nabla \cdot \bp$,    and applying integration by parts, using H\"older's inequality along with  \eqref{1.2} 
 we calculate 
 \begin{align}
&  \intl_{B(4r) \setminus B(r)} |u|^3 \psi^3   dx   
\cr
& \quad = \intl_{B(4r) \setminus B(r)} \partial _j (\Phi_{ ij}- (\Phi_{ ij})_{ B(4r)}) u _i |u| \psi^3   dx  
\cr
&\quad  = -\intl_{B(4r) \setminus B(r)}  (\Phi_{ ij} - (\Phi_{ ij})_{ B(4r)}) \partial _j (u _i |u|) \psi^3   dx  -  
\intl_{B(4r) \setminus B(r)}   (\Phi_{ ij}- (\Phi_{ ij})_{ B(4r)}) (u _i |u|)\partial _j \psi^3   dx 
\cr
& \quad \le  c \bigg(\intl_{B(4r) }   |\bp- \bp_{ B(4r)}|^6  dx\bigg)^{ \frac{1}{6}}
\bigg(\intl_{B(4r) \setminus B(r)}  |u|^3 \psi^3   dx \bigg) ^{ \frac{1}{3}}\bigg(\intl_{B(4r) \setminus B(r)}  |\nabla u|^2    dx \bigg) ^{ \frac{1}{2}}
\cr
& \quad+   c r^{ -1}\bigg(\intl_{B(4r) }   |\bp- \bp_{ B(4r)}|^6  dx\bigg)^{ \frac{1}{6}}
\bigg(\intl_{B(4r) \setminus B(r)}  |u|^3 \psi^3   dx \bigg) ^{ \frac{1}{3}}\bigg(\intl_{B(4r) \setminus B(r)}  | u|^2    dx \bigg) ^{ \frac{1}{2}}
\cr
&\quad  \le c r^{ \frac{2}{3}} \bigg(\intl_{B(4r) \setminus B(r)}  |u|^3 \psi^3   dx \bigg) ^{ \frac{1}{3}}\bigg(\intl_{B(4r) \setminus B(r)}  |\nabla u|^2    dx \bigg) ^{ \frac{1}{2}} 
\cr
& \qquad \quad +   c r^{ - \frac{1}{3}}\bigg(\intl_{B(4r) \setminus B(r)}  |u|^3 \psi^3   dx \bigg) ^{ \frac{1}{3}}\bigg(\intl_{B(4r) \setminus B(r)}  | u|^2    dx \bigg) ^{ \frac{1}{2}}. 
 \label{2.9}
\end{align}

First, we note that  \eqref{2.0e} with $ R=4r$ and $ \varphi = \psi $ implies 
\begin{align}
|u_{ B(4r),  \psi }|  &\le   c r^{  - \frac{5}{6}}. 
 \label{2.8}
\end{align}
By the triangular inequality and  \eqref{2.0d} with $ E= B(4r) \setminus B(r)$ we have
\begin{align*}
\bigg(\intl_{B(4r) \setminus B(r)}  | u|^2    dx \bigg) ^{ \frac{1}{2}} &\le \bigg(\intl_{B(4r) \setminus B(r)}  |u- u_{ B(4r), \psi }|^2 dx  \bigg) ^{ \frac{1}{2}} 
+ \bigg(\intl_{B(4r) \setminus B(r)} |u_{ B(4r), \varphi  } |^2 dx  \bigg) ^{ \frac{1}{2}}
\cr
& \le c \bigg(\intl_{B(4r) \setminus B(r)}  |u- u_{B(4r) \setminus B(r) }|^2 dx  \bigg) ^{ \frac{1}{2}} 
+ c r^{ \frac{3}{2}} |u_{ B(4r),  \psi }|. 
\end{align*}
Using the Poincar\'{e} inequality and \eqref{2.8},  we find 
\be\label{2.8c}
\bigg(\intl_{B(4r) \setminus B(r)}  | u|^2    dx \bigg) ^{ \frac{1}{2}}  \le c r\bigg(\intl_{B(4r) \setminus B(r)}  |\nabla u|^2    dx \bigg) ^{ \frac{1}{2}} + cr^{ \frac{2}{3}}.
\ee
 Inserting \eqref{2.8c}  into the last term of  \eqref{2.9},   and dividing the result by $ \bigg(\intl_{B(4r) \setminus B(r)}  |u|^3 \psi^3   dx \bigg) ^{ \frac{1}{3}}$, we find 
\begin{align}\label{2.8d}
  \intl_{B(4r) \setminus B(r)} |u|^3 \psi^3   dx   
 \le cr \bigg(\intl_{B(4r) \setminus B(r)}  |\nabla u|^2    dx \bigg) ^{ \frac{3}{4}} + c r^{ \frac{1}{2}}. 
\end{align}
Repeating the estimates from \eqref{2.9}  to \eqref{2.8d} for $B$,  we obtain
\begin{align*}
  \intl_{B(4r) \setminus B(r)} |B|^3 \psi^3   dx   
 \le cr \bigg(\intl_{B(4r) \setminus B(r)}  |\nabla B|^2    dx \bigg) ^{ \frac{3}{4}} + c r^{ \frac{1}{2}}. 
\end{align*}

Thus,    observing \eqref{2.4},  we obtain the claim  \eqref{2.6}.   \\

Let $ 1 < r< +\infty$ be arbitrarily chosen.   By  $ \zeta \in  C^\infty(\R^{n} )$ we denote  a cut off function, which is radially non-increasing   with $ \zeta =1$ on $B(2r)$ and $\zeta=0$ on $ \Bbb R^3\setminus B(3r)$ such that  
$ |\nabla \zeta | \le c r^{ -1}$ and $ |D ^2 \zeta | \le c r^{ -2}$.  We multiply  \eqref{21} by 
$ u \zeta $, and  integrate it over $ B(3r)$ ,   and similarly  multiply \eqref{22} by $B\zeta$ and integrate it over $B(3r)$, and adding them together, we find after integration by parts. 
\begin{align}
& \intl_{B(3r) } |\nabla u|^2 dx + \intl_{B(3r) } |\nabla B|^2 dx  \le  \frac{1}{2}\intl_{B(3r) \setminus B(2r)} |u|^2 \Delta \zeta ^2 dx \cr
 &\qquad + \frac{1}{2}\intl_{B(3r) \setminus B(2r)} |u|^2 u \cdot \nabla \zeta ^2  dx+ \intl_{B(3r) \setminus B(2r)} (p - p_{ B(\overline{R})}) u \cdot \nabla \zeta ^2 dx \cr
 &\qquad - \intl_{B(3r) \setminus B(2r)}  B\cdot u  ( B\cdot \nabla ) \zeta^2 dx +  \frac{1}{2}\intl_{B(3r) \setminus B(2r)} |B|^2 \Delta \zeta ^2 dx \cr
 &\qquad  + \frac{1}{2}\intl_{B(3r) \setminus B(2r)} |B|^2 u \cdot \nabla \zeta ^2  dx    + \intl_{B(3r) \setminus B(2r)}  u\cdot B   ( B\cdot \nabla ) \zeta^2 dx \cr
&\qquad -\alpha \intl_{B(3r) \setminus B(2r)} ((\nabla \times B) \times B) \cdot B\times \nabla \zeta^2 dx \cr
      &\le cr^{ -2} \intl_{B(3r) \setminus B(2r)} |u|^2  dx+  cr^{ -1} \intl_{B(3r) \setminus B(2r)} |u|^3dx 
       + c  r^{ -1} \intl_{B(3r) \setminus B(2r)} |p - p_{ B(\overline{R}) }| |u| dx\cr
       & + c  r^{ -1} \intl_{B(3r) \setminus B(2r)} |B|^2 |u| dx + c r^{ -2} \intl_{B(3r) \setminus B(2r)} |B|^2  dx
    +\alpha  c  r^{ -1}    \intl_{B(3r) \setminus B(2r)} |B|^2 |\nabla B| dx\cr
  &= J_1+ \cdots +J_6. 
 \label{2.1d}
\end{align}
Using    \eqref{2.5},  we immediately get 
\[
J_1 +J_5\le c r^{ -1} \bigg(\intl_{B(3r)} |u|^3 dx \bigg)^{ \frac{2}{3}} +   c r^{ -1} \bigg(\intl_{B(3r)} |B|^3 dx \bigg)^{ \frac{2}{3}} \le c r^{ - \frac{1}{3}}= o(1)\quad  \text{as}\quad r \rightarrow +\infty. 
\]
From  \eqref{2.6} it follows that $ J_2= o(1)$ as $ r \rightarrow +\infty$. 
For $J_3$ we observe 
\begin{align*}
J_3& \le c r^{ -1}\bigg( \intl_{B(3r)} |p - p_{ B(3r)}|^{ \frac{3}{2}} dx  \bigg)^{ \frac{2}{3}}\bigg( \intl_{B(3r)\setminus B(2r) } |u|^{3} dx  \bigg)^{ \frac{1}{3}} \cr
&\le  c\bigg( r^{ -1} \intl_{B(3r)} |p - p_{ B(3r)}|^{ \frac{3}{2}} dx  \bigg)^{ \frac{2}{3}} \bigg(r^{-1} \intl_{B(3r)\setminus B(2r) } |u|^{3} dx  \bigg)^{ \frac{1}{3}}\cr
&\le c \bigg( r^{ -1} \intl_{B(3r)} |p - p_{ B(3r)}|^{ \frac{3}{2}} dx  \bigg)^{ \frac{2}{3}} o(1). 
\end{align*}
Using the estimate  \eqref{2.21} for $ 3r$ in place of $ \overline{R}$, we obtain 
\begin{align*}
\hspace*{-1cm}   r^{ -1} \intl_{B(3r)} |p - p_{ B(3r)}|^{ \frac{3}{2}} dx  &\le c r^{ -1} \intl_{B(3r)} |\nabla u|^{ \frac{3}{2}} dx + 
   c r^{ -1}  \intl_{B(3r)} |u|^3 dx +c r^{ -1}  \intl_{B(3r)} |B|^3 dx.\\[0.2cm]
&\le c r^{ - \frac{1}{4}} \bigg(\intl_{B(3r)}  |\nabla u|^2 dx \bigg)^{ \frac{3}{4}}  +  c r^{ -1}  \intl_{B(3r)} |u|^3 dx+c r^{ -1}  \intl_{B(3r)} |B|^3 dx. 
\end{align*}
By virtue of  \eqref{2.4} and  \eqref{2.5} the right-hand side of the above inequality is bounded for  $ r \ge 1 $.  This shows that $ J_3 =o(1)$ as $ r 
\rightarrow +\infty$. 

For $J_4$  by the H\"{o}lder inequality and \eqref{2.6} we have
\begin{align*}
J_4 & \le c r^{ -1}\bigg( \intl_{B(3r)\setminus B(2r)} |B|^3 dx  \bigg)^{ \frac{2}{3}}\bigg( \intl_{B(3r)\setminus B(2r) } |u|^{3} dx  \bigg)^{ \frac{1}{3}} \cr
&\le  c\bigg( r^{ -1} \intl_{ B(3r)\setminus B(2r) }|B|^3 dx  \bigg)^{ \frac{2}{3}} \bigg(r^{-1} \intl_{B(3r)\setminus B(2r) } |u|^{3} dx  \bigg)^{ \frac{1}{3}}=o(1)
\end{align*}
as $r \to +\infty$.
Finally for $J_6$,    by the H\"{o}lder and the Sobolev-Poincar\'{e} inequalities  together with \eqref{2.0e},  we have
\begin{align*}
J_6& \le \alpha  c  r^{ -1}   \bigg( \intl_{B(3r) \setminus B(2r)} |B|^{4}  dx \bigg)^{\frac12}  \bigg( \intl_{B(3r) \setminus B(2r)}  |\nabla B| ^2 dx\bigg)^{\frac12}\cr
& \le \alpha  c  r^{ -\frac12}     \bigg( \intl_{B(3r) \setminus B(2r)} |B|^{6}  dx \bigg)^{\frac13}  \bigg( \intl_{B(3r) \setminus B(2r)}  |\nabla B| ^2 dx\bigg)^{\frac12}\cr
&  \le \alpha  c  r^{ -\frac12}     \bigg( \intl_{B(3r) \setminus B(2r)} |B- B _{ B(3r) \setminus B(2r) , \zeta}  |^{6}  dx \bigg)^{\frac13}  \bigg( \intl_{B(3r) \setminus B(2r)}  |\nabla B| ^2 dx\bigg)^{\frac12}\cr
 &\qquad +c \alpha  r^{ \frac12}    | B _{ B(3r) \setminus B(2r) , \zeta}|^2 \bigg( \intl_{B(3r) \setminus B(2r)}  |\nabla B| ^2 dx\bigg)^{\frac12} 
 \end{align*}
\begin{align*}
& \le \alpha  c  r^{ -\frac12}   \bigg( \intl_{B(3r) \setminus B(2r)} |\nabla B|^{2}  dx \bigg)^ {\frac32}   + c \alpha  r^{-\frac76} \bigg( \intl_{B(3r) \setminus B(2r)}  |\nabla B| ^2 dx\bigg)^{\frac12}\cr
&=o(1)
\end{align*}
as $r\to +\infty$.
Inserting  the above estimates of $J_1, \cdots,J_6 $ into the right-hand side of  \eqref{2.1d}, we deduce that 
\[
\intl_{B(r)} |\nabla u|^2 dx + \intl_{B(r)} |\nabla B|^2 dx  = o(1)\quad   \text{as}\quad  r \rightarrow +\infty. 
\]
Accordingly, $ u \equiv const$ and $ B \equiv const$ and by means of  \eqref{2.5} it follows that $u=B= 0$. 
 \hfill \Beweisende

\hspace{0.5cm}

$$\mbox{\bf Acknowledgements}$$
Chae was partially supported by NRF grant 2016R1A2B3011647, while Wolf has been supported 
supported by NRF grant 2017R1E1A1A01074536.
The authors declare that they have no conflict of interest.
\bibliographystyle{siam}

\end{document}